\documentclass[12pt]{article}
\usepackage{rotating}
\usepackage{graphicx}
\usepackage{latexsym}
\usepackage{amsfonts,amsthm,bm}

\newcommand{\TC}{\ensuremath{\mathcal{T}}}
\newcommand{\RTC}{\ensuremath{\mathcal{R}}}

\newcommand{\peq}{\mbox{$\,:\stackrel{\scriptscriptstyle \,+}{=}\,$}}

\newcommand{\rmd}{\mbox{\rm d}}
\newcommand{\rmi}{\mbox{\rm i}}

\begin{document}

\title{Fuchsian differential equation for the perimeter generating function of
three-choice polygons}

\author{Anthony J. Guttmann and Iwan Jensen \\
\small ARC Centre of Excellence for Mathematics and Statistics of Complex Systems, \\
\small Department of Mathematics and Statistics, \\
\small The University of Melbourne, Victoria 3010, Australia}

\date{\today}

\maketitle

\begin{abstract}
Using a simple transfer matrix approach we have derived very long series expansions 
for the perimeter generating function of three-choice  polygons.
We find that all the terms in the generating function can be reproduced from a
linear Fuchsian differential equation of order 8. We perform an analysis of the
properties of the differential equation.
\end{abstract}

\section{Introduction}

A well-known long standing problem in combinatorics and statistical mechanics is
to find the generating function for self-avoiding polygons (or walks) on a two-dimensional
lattice, enumerated by perimeter. Recently, we have gained a greater understanding 
of the difficulty of this problem, as Rechnitzer \cite{AR03a} has {\em proved} that the 
(anisotropic) generating function for square lattice self-avoiding polygons is not 
differentiably finite \cite{RPS80a}, as had been previously {\em conjectured}, on numerical 
grounds \cite{GC01}, but not proved. That is to say, it cannot be expressed as the solution 
of an ordinary  differential equation with polynomial coefficients. There are many 
simplifications of this problem that are solvable \cite{BM96a}, but all the simpler models 
impose an effective directedness or other constraint that reduces the problem, in essence, 
to a one-dimensional problem.

One model, that of so-called {\em three-choice polygons}, has remained unsolved despite the
knowledge that its solution must be D-finite. In this paper we report on recent numerical work
resulting in an exact differential equation apparently satisfied by the perimeter
generating function of three-choice polygons. While our results do not constitute a rigorous
mathematical proof the numerical evidence is compelling.

Three-choice self-avoiding walks on the square lattice, ${\mathbb Z}^2$,
were introduced by Manna \cite{Manna84} and can be
defined as follows: Starting from the origin one can step in any direction; after a step
upward or downward one can head in any direction (except backward); after a step
to the left one can only step forward or head downward, and similarly after a step to 
the right one can continue forward or turn upward. Alternatively put, one cannot make a right-hand
turn after a horizontal step. Whittington \cite{SGW85a} showed that
the growth constant for three-choice walks is exactly 2, so that if $w_n$ denotes
the number of such walks of $n$ steps  on an infinite lattice, equivalent up to a translation,
then $w_n \sim 2^{n + {\rm o}( n)}.$  It is perhaps surprising that the best known result for the
sub-dominant term is $2^{{\rm o}(n)}$ but attempts to improve on this have not been successful.
Even numerically, there is no firmly based conjecture for the sub-dominant term, unlike for
ordinary self-avoiding walks, for which the sub-dominant term is widely believed to be ${\rm O}(n^g).$

\begin{figure}
\begin{center}
\includegraphics{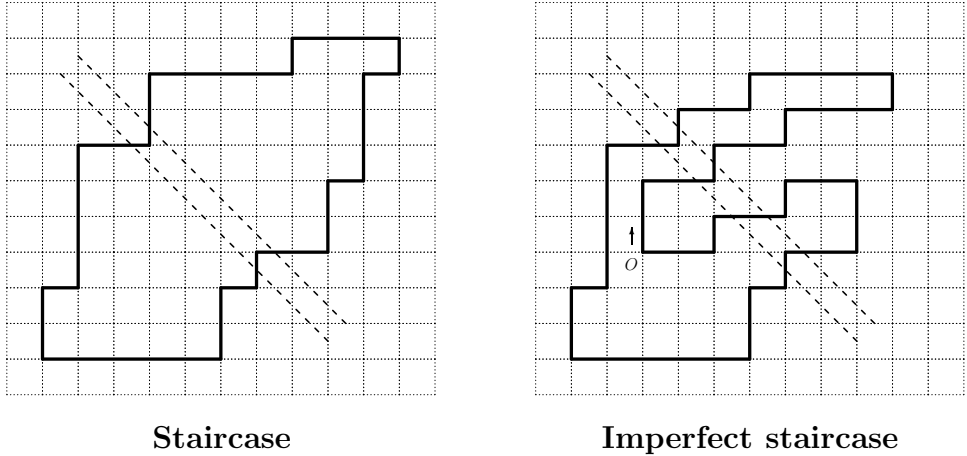}
\end{center}
\caption{\label{fig:poly} Examples of the two types of three-choice polygons. In the right
panel we indicate the origin ($O$) and the direction of the first step (note that rotation by
180 degrees also leads to a valid three-choice polygon). 
}
\end{figure}

As usual one can define a polygon version of the walk model by requiring the walk to return to the origin. 
So a three-choice polygon \cite{GPO93} is simply a three-choice self-avoiding walk which returns to the 
origin, but has no other self-intersections. There are two distinct classes of three-choice polygons. The 
three-choice rule either leads to staircase polygons or {\em imperfect staircase polygons} \cite{CGD97} 
as illustrated in figure~\ref{fig:poly}. In the case of staircase polygons any vertex on the perimeter can 
act as the origin of the three-choice walk (which then proceeds counter-clockwise), while for imperfect 
staircase polygons there is only one possible origin but the polygon could be rotated by 180 degrees. If we 
denote by $t_n$ the number of three-choice polygons with perimeter $2n$ then, $t_n = 2nc_n + 2p_n$,  where 
$c_n$ is the number of staircase polygons with perimeter $2n$, and $p_n$ is the number of imperfect 
staircase polygons of perimeter $2n$. Note that $t_n$, $p_n$ and $c_n$ all grow like $4^n$ and
in particular we recall the well-known result that $c_{n+1}=C_n=\frac{1}{n+1} {2n \choose n}$ are given by 
the Catalan numbers $C_n$.

In this paper we report on recent work which has led to an exact Fuchsian \cite{Ince} linear
differential equation of order 8 apparently satisfied by the perimeter generating function, 
$\TC(x) = \sum_{n\geq 0} t_nx^n$, for three-choice polygons (that is $\TC (x)$ is conjectured to be
one of the solutions of the ODE, expanded around the origin). The first few terms in the
generating function are 
$$\TC(x) = 4x^2 + 12x^3+42x^4+152x^5+562x^6+\cdots.$$ The generating function for the coefficients
$p_n$ is no simpler.

If we distinguish between steps in the $x$ and $y$ direction, and let $t_{m,n}$ denote
the number of three-choice polygons with $2m$ horizontal steps and $2n$ vertical steps,
then the anisotropic generating function for $\TC$ can be written
$$\TC (x,y) = \sum_{m,n} t_{m,n}x^m y^n = \sum_n H_n(x)y^n,$$ where $H_n(x) = \frac{R_n(x)}{S_n(x)}$
is the (rational \cite{RPSECv2}) generating function for three-choice polygons with $2n$ vertical steps.
In earlier, unpublished, numerical work, we found that, for imperfect staircase polygons, the 
denominators were:
$$S_n(x) = (1 - x)^{2n - 1}(1 + x)^{(2n-7)_+} \qquad n \quad {\rm even,}$$ and 
$$S_n(x) = (1 - x)^{2n - 1}(1 + x)^{(2n-8)_+} \qquad n \quad {\rm odd}.$$
This was subsequently proved by Bousquet-M\'elou \cite{BMb}. 
Further, Bousquet-M\'elou showed that the numerators satisfied:
$$R_n(-1)=\frac{-12(4m)!}{m!(m + 1)!(m + 2)!(m + 3)!} \qquad n=2m+4,$$ and
$$R_n(-1)=\frac{-96(4m+1)!}{m!(m + 1)!(m + 2)!(m + 4)!} \qquad n=2m+5.$$ 
Unfortunately, we still do not have enough information to identify the numerators, though
we observe that they are of degree $3n - 7$ for $n \ge 4$ and $n$ even, and of degree
$3n-8$ for $n \ge 5$ and $n$ odd.

\begin{figure}
\begin{center}
\includegraphics[scale=0.5]{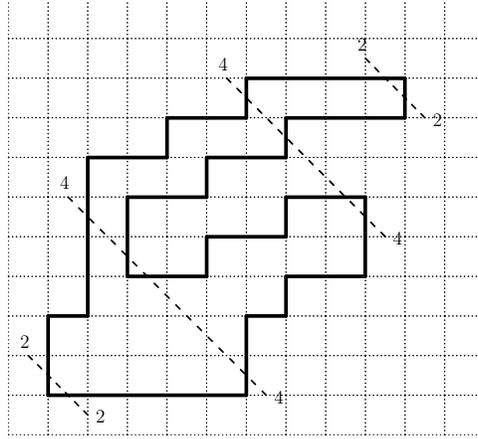}
\end{center}
\caption{ \label{fig:gv} Showing the decomposition of an imperfect staircase polygon
into a sequence of 2-4-2 non-intersecting walkers, each expressible as a Gessel-Viennot determinant}
\end{figure}

It is also possible to express the generating function $\TC(x)$ as a five-fold sum, with one 
constraint, \cite{BMb} of $4 \times 4$ Gessel-Viennot determinants \cite{GV89}. This is clear 
from figure \ref{fig:gv}, where the enumeration of the lattice paths between the dotted lines is just
the classical problem of 4 non-intersecting walkers, and these must be joined to two non-intersecting 
walkers to the left, and to two non-intersecting walkers to the right. Then one must sum over 
different possible geometries. However the fact that the generating function is so expressible 
implies that it is differentiably finite \cite{L89}. 

In the following sections we report on our work leading to an ODE for the perimeter
generating function of three-choice polygons. 
We started by generating the counts for three-choice polygons up to half-perimeter 260.
Using numerical experimentation we then found what we believe to be the underlying
ODE. This calculation required the use of the first 206 coefficients with 
the resulting ODE then correctly predicting the next 54 coefficients.
While the possibility that this ODE is not the correct one is extraordinarily
small, our result does of course not constitute a proof.
Unfortunately we cannot usefully bound the size of the underlying ODE, otherwise we could
use the knowledge of D-finiteness to provide a proof of our results. That is to say, any bounds
that follow from closure theorems \cite{L89} are too large to be useful.

\section{Computer enumeration \label{sec:enum}}

The algorithm we use to count the number of imperfect polygons is a slightly modified 
version of the algorithm of Conway {\it et al.} \cite{CGD97}. Before proceeding to the full 
problem it is useful to briefly outline the transfer matrix algorithm for enumerating 
staircase polygons. Recall that a staircase polygon consists of two directed walks starting 
at the origin, moving only to the right and up, and terminating once the walks join at a 
vertex. If we look at a diagonal line $x+y=k+1/2$ then for any integer $k$ this line will 
intersect a polygon at 0 (miss the polygon) or 2 edges (intersect the polygon), 
see figure~\ref{fig:poly}. We start with $k=0$ such that the line intersects the first two 
edges of the staircase polygon. We then move the line upward (increase $k$ by 1) and as we 
do this we add an edge to each walk. There are only four new configurations corresponding to 
the four possible steps. We need only keep track of the gap between the two walks, where the 
gap is the minimal number of iterations required in order to join the two walks. As we move the 
line the gap is either increased by a unit (the upper walk moves up and the lower walk moves right), 
decreased by a unit  (the upper walk moves right and the lower walk moves up) or remains constant 
in two possible ways (both walks move up or right). Let $C(i,k)$ be the number of configurations
with a gap of $i$ at step $k$. We then have the following very simple algorithm:

Set $C(1,0) = x$ (where $x$ is a variable conjugate to the half-perimeter of the polygon). 
Run through all possible gaps $i=1,\ldots, k+1$ and do the
following updates: $C(i+1,k+1)\peq xC(i,k)$, $C(i-1,k+1)\peq xC(i,k)$ and
$C(i,k+1)\peq 2xC(i,k)$. 
Here $a\peq b$ is short-hand for assign to $a$ the value $a+b$.

Formally we can view the transformation from the set of states $C(i,k)$ to $C(j,k+1)$ as a 
matrix multiplication (hence our use of the nomenclature {\em transfer matrix algorithm}) 
with $k$ counting the number of iterations of the transfer matrix algorithm. However, as can 
be readily seen from the algorithm the transfer matrix is extremely sparse and there is no 
reason to list it explicitly (it is given {\em implicitly} by the updating rules).

The term $C(0,k)$ is the number of staircase polygons of half-perimeter $k+1$.
Note that the use of the variable $x$ is somewhat 
superfluous in the case of staircase polygons since the generating function
at iteration $k$ is just  $x^{k+1}C(i,k)$, where $C(i,k)$ is the number of
configuration with gap $i$ after $k$ iterations. But it is included here for
reasons of generality and in the case of imperfect staircase polygons the 
generating function will be a (non-trivial) polynomial in $x$.
Naturally we need not actually keep all the entries $C(i,k)$ since only the
current and subsequent values are needed for the calculation so we can replace $C(i,k)$ 
with $C(i,k \bmod 2)$. We just have to initially set to zero all entries in the next 
step and keep a running total $c(k)$ of the number of staircase polygons. 

Imperfect staircase polygons start out as ordinary staircase polygons (see figure~\ref{fig:gv})
Then at some vertex two additional directed walks (sharing the same starting point) 
are inserted between the two original walks (at the first dashed line marked with a `4' in 
figure~\ref{fig:gv}). The diagonal line will thus intersect these polygon configurations at 
four edges. Imperfect staircase polygons are created by connecting the first two walks and 
the last two walks (as illustrated at the last dashed line marked with a `4' in figure~\ref{fig:gv}). 
With four walks we need to retain three pieces of information, namely, the three gaps $l$, $m$, 
and $n$ between consecutive walks. Each existing configuration can produce 16 new configurations
as each walk is extended by a step either up or to the right. The resulting updating 
is easily worked out \cite{CGD97}. Let $G(l,m,n)$ be the generating function (a polynomial
in the variable $x$) for partially completed polygons at a given diagonal $k$. As we proceed 
to the next diagonal $k+1$ we add $x^2G(l,m,n)$ (the factor $x^2$ arise because we extend all 
walks by a step) to $G(l,m,n)$ (twice), $G(l+1,m,n)$, $G(l,m+1,n)$,  $G(l,m,n+1)$, $G(l-1,m,n)$,  
$G(l,m-1,n)$, $G(l,m,n-1)$, $G(l+1,m-1,n)$, $G(l+1,m,n-1)$, $G(l-1,m+1,n)$, $G(l,m+1,n-1)$, 
$G(l-1,m,n+1)$, $G(l,m-1,n+1)$, $G(l-1,m+1,n-1)$  and $G(l+1,m-1,n+1)$. 
Any update resulting in $G(l,0,n)$  has to be rejected because it corresponds to a configuration
in which we have joined the two middle walks in and this can never
lead to an imperfect staircase polygon. Obviously once any two walks have been
connected the remaining walks follow the usual staircase polygon updating rules.

The configurations with two walks already connected can also be encoded by the $G$ functions.
We simply let $G(l,0,0)$ be the generating function for partial polygons with two walks already 
connected (note that if the boundary line intersects four edges $m>0$). So in the updating of 
imperfect staircase polygons we can set $G(0,m,n)$ (we connect the two lower-most walks) to 
$G(n,0,0)$ (this case is illustrated at the last dashed line marked with a `4' in figure~\ref{fig:gv}).
Likewise we can set $G(l,m,0)$ (we connect the two uppermost walks) to $G(l,0,0)$.
The condition for the formation of a valid polygon is $l=n=0$ (note that we can't demand 
$m=0$ as well, since we could connect {\em both} the lower- and uppermost walks simultaneously).

The `creation' of a configuration with three gaps, alternatively, one in which a diagonal line 
intersects four edges of an imperfect staircase polygon, which we refer to as a $G$-type configuration,
is also very simple (see the first dashed line marked with a `4' in figure~\ref{fig:gv}). We start 
with a staircase type configuration $C(i)$ and from this we can create four $G$-type configurations 
by assigning the value $x^2C(i)$ to $G(j,1,i-j)$, $G(j-1,1,i-j)$, $G(j,1,i-j-1)$ and $G(j-1,1,i-j-1)$, 
where $1 \leq j \leq i-1$ (the factor $x^2$ arises because we extend the outer walks by a step {\em and}
insert two new walks each containing a single step). 

The algorithm outlined above is already very efficient, but it can be further enhanced
by the following simple observation. If we wish to calculate the number of polygons
up to a given maximal half-perimeter length $N$, we need not consider all possible gaps
since some configurations can only lead to polygons of a size exceeding $N$. First
of all since gaps only increase or decrease by one at each iteration we need never
consider configurations with gaps exceeding $N/2$. Furthermore,  any $G$-type configuration
with $m>0$ must have half-perimeter at least $k+m$. Here we get the contribution $k$ from the 
outermost walks ($k$ is the number of forward steps or iterations taken) and the contribution $m$
from the innermost walks (a gap $m$ requires at least $m$ steps). In order to produce an 
imperfect staircase polygon we have to add at least $l+n$ additional steps (we have to join both
the two upper-most and two lower-most walks), so if $M=k+l+m+n>N$ we can discard this
configuration. Not only can we thus discard some configurations when $M>N$ but we can 
also further decrease the memory use since rather than storing $N$ terms per configuration 
we only need to store $N-M$ terms. 

We calculated the number of imperfect staircase polygons up to perimeter 520.  
The integer coefficients become very large so the calculation was performed using 
modular arithmetic \cite{KnuthACPv2}. This involves performing the calculation modulo 
various prime numbers $p_i$ and then reconstructing the full integer coefficients 
at the end. We used primes of the form $p_i=2^{30}-r_i$, where $r_i$ are small positive
integers, less than $1000,$ chosen so that $p_i$ is prime, and $p_i \ne p_j$ unless $i = j.$
18 primes were needed to represent the coefficients correctly. The calculation for each 
prime used about 250Mb of memory and about 18 minutes of CPU time on a 
2.8 GHz Xeon processor. Naturally we could have carried the calculation much further
but as we shall demonstrate in the next section this more than sufficed to 
identify an exact differential equation satisfied by $\TC (x)$.

\section{The Fuchsian differential equation \label{sec:fde}}

In recent papers Zenine {\it et al.} \cite{ZBHM04a,ZBHM05a,ZBHM05b} obtained the
linear differential equations whose solutions give the 3- and 4-particle contributions
$\chi^{(3)}$ and $\chi^{(4)}$ to the Ising model susceptibility. In this paper we use their
method to find a linear differential equation which has as a solution the generating 
function $\TC (x)$ for three-choice polygons. We briefly outline the method here. Let us
assume we have a function $F(x)$ with a singularity at $x=x_c=1/\mu$. Starting from a (long) series 
expansion for the function $F(x)$ we look for a linear differential equation of order $m$ 
of the form

\begin{equation}\label{eq:de}
\sum_{k=0}^m P_k(x) \frac{\rmd^k}{\rmd x^k}F(x) = 0,
\end{equation}
such that $F(x)$ is a solution to this homogeneous linear differential equation, where
the $P_k(x)$ are polynomials. In order to make it as simple as possible we start
by searching for a Fuchsian \cite{Ince} equation. Such equations have only regular singular points. There
are several reasons for searching for a Fuchsian equation, rather than a more general D-finite
equation. Computationally the Fuchsian assumption simplifies the search for a  solution. 
One may also argue, less precisely, that for ``sensible''  combinatorial models one
would expect Fuchsian equations, as irregular singular points are characterized by explosive,
super-exponential behaviour. Such behaviour is not normally characteristic of combinatorial
problems. (The point at infinity may be an exception to this somewhat imprecise observation).
One may also ask the question whether most of the problems in combinatorics with D-finite
solutions have Fuchsian solutions? While we have not made an exhaustive study, we know of
no counter-example to this suggestion.  

From the general theory of Fuchsian \cite{Ince} equations it follows  that the degree of 
$P_k(x)$ is at most $n-m+k$ where $n$ is the degree of $P_m(x)$. To simplify matters (reduce 
the order of the unknown polynomials) it is often advantageous to explicitly assume that the 
origin and $x=x_c$ are regular singular points and to set $P_k(x)=Q_k(x)S(x)^k$, where $S(x)=xR(x)$ 
and $R(x)$ is a polynomial of minimal degree having $x_c$ as a root (in our case we have $R(x)=1-4x$). 
$S(x)$ could be generalised to include more regular singular 
points if some were known from other methods of analysis, but we have not found this to be
advantageous. Thus when searching for a solution of Fuchsian type there are 
only two parameters: namely the order $m$ of the differential equation and the degree $q_m$ 
of the polynomial $Q_m(x)$. Let $\rho$ be the degree of $S(x)$ (2 in our case), then for given 
$m$ and $q_m$ there are $L=(m+1)(q_m+1)+\rho m(m+1)/2-1$ unknown coefficients, where we  have 
assumed without loss of generality that the leading order coefficient in $P_m(x)=Q_m(x)S(x)^m$ is 1. 
We can then search systematically for solutions by varying $m$ and $q_m$. In this
way we first found a solution with $m=10$ and $q_m=12$, which required the determination
of $L=206$ unknown coefficients. We have 260 terms in the half-perimeter series and
thus have more than 50 additional terms with which to check the correctness of our solution.
Having found this conjectured solution we then turned the ODE into a recurrence relation and used this
to generate more series terms in order to search for a lower order Fuchsian equation.
The lowest order equation we found was eighth order and with $q_m=30$, which requires the determination
of $L=321$ unknown coefficients. Thus from our original 260 term series we could not have found
this $8^{th}$ order solution since we did not have enough terms to determine all the unknown coefficients
in the ODE. This raises the question as to whether perhaps there is an ODE of lower order than 8 that
generates the coefficients? The short answer to this is no. Further study \cite{GJ06} of our differential 
operator revealed that it can be factorised. In fact we found a factorization into three first-order
linear operators, a second order and a third order. The generating function is a solution of
the $8^{th}$ order operator, not of any of the smaller factors.

So the (half)-perimeter generating function $\TC (x)$ for three-choice polygons
is conjectured to be a solution of the linear differential equation of order 8

\begin{equation}\label{eq:TCfde}
\sum_{k=0}^8 P_k(x) \frac{\rmd^k}{\rmd x^k}F(x) = 0
\end{equation}

with

\begin{eqnarray}
& & P_8(x)=x^3(1-4x)^4(1 + 4x)(1 + 4x^2)(1 + x + 7x^2)Q_8(x),   \nonumber \\
& & P_7(x)=x^2(1-4x)^3 Q_7(x), \;\;\;\;\;\; P_6(x)=2x (1-4x)^2 Q_6(x),   \nonumber \\
& & P_5(x)=6(1-4x) Q_5(x), \,\,\;\;\; P_4(x)=24 Q_4(x),  \label{eq:TCpol} \\
& & P_3(x)=24 Q_3(x), \,\;\; P_2(x)=144x(1-2x) Q_2(x),    \nonumber \\
& & P_1(x)=144(1-4x) Q_1(x), \; P_0(x)=576 Q_0(x),   \nonumber
\end{eqnarray}
where $Q_8(x)$, $Q_7(x)$, $\ldots$, $Q_0(x)$, are polynomials of degree 25, 31, 32, 33, 33,
32, 29, 29, and 29, respectively. The polynomials are listed in Appendix~A (note that
the polynomials do not factorise).

The singular points of the differential equation are given by the roots of $P_8(x)$.
One can easily check that all the singularities (including $x=\infty$) are
{\em regular singular points} so equation (\ref{eq:TCfde}) is indeed of the Fuchsian type.
It is thus possible, using the method of Frobenius, to obtain from the indicial equation
the critical exponents at the singular points. These are listed in Table~\ref{tab:TCexp}.

\begin{table}
\caption{\label{tab:TCexp}
Critical exponents for the regular singular points of the Fuchsian differential
equation satisfied by $\TC (x)$.}
\begin{center}
\begin{tabular}{ll}
\hline \hline
Singularity & Exponents \\
\hline
$x=0$ & $-1, \, 0, \, 0, \, 0, \, 1, \, 2, \, 3, \, 4$ \\
$x=1/4$ & $-1/2, \, -1/2, \, 0, \, 1/2, \, 1, \, 3/2, \, 2, \, 3$ \\
$x=-1/4$ & $0, \, 1, \, 2, \, 3, \, 4, \, 5, \, 6, \, 13/2$ \\
$x=\pm\, \rmi/2$ & $0, \, 1, \, 2, \, 3, \, 4, \, 5, \, 6, \, 13/2$ \\
$1+x+7x^2=0$ & $0, \, 1, \, 2, \, 2, \, 3, \, 4, \, 5, \, 6$ \\
$x=\infty$ & $-2, \, -3/2, \, -1, \, -1, \, -1/2, \, 1/2, \, 3/2, \, 5/2$ \\
$Q_8(x)=0$ & $0, \, 1, \, 2, \, 3, \, 4, \, 5, \, 6, \, 8$ \\
\hline \hline
\end{tabular}
\end{center}
\end{table}

For equations of the Fuchsian type the critical exponents satisfy a simple Fuchsian summation 
relation, which we now take the opportunity to confirm in our case. 
Let $x_1,\, x_2,\ldots x_n,\, x_{n+1}=\infty$ be the regular singular points of
a Fuchsian type equation of order $m$ and $\alpha_{j,1}, \ldots \alpha_{j,m}$ $(j=1,\ldots n+1)$
the $m$ exponents determined from the roots of the indicial equation for each regular singular
point, $x_j$, then the following Fuchsian relation holds:
\begin{equation}
\sum_{j=1}^{n+1}\sum_{k=1}^m \alpha_{j,k}= \frac{(n-1)m(m-1)}{2}.
\end{equation}
In this case the number of regular singular points is $m+1=33$, namely the 25 roots
of $Q_8(x)$, the two roots of $1+x+7x^2$, $x=\pm\, \rmi/2$, $x= \pm 1/4$,
$x=0$ and $x=\infty$. It is easy to verify that the Fuchsian relation is satisfied
with $m=8$, $n=32$, and all the exponents $\alpha_{j,k}$ summing to 868, which is a useful
check on our results.

We shall now consider the local solutions of the differential equation around each singularity. 
Recall that in general it is known \cite{ForsythV4,Ince} that if the indicial equation yields $k$ 
critical exponents which differ by an integer, then the local solutions {\em may} 
contain logarithmic terms up to $\log^{k-1}$. However, for the Fuchsian equation (\ref{eq:TCfde})
{\em only} multiple roots of the indicial equation give rise to logarithmic terms in the local 
solution around a given singularity, so that a root of multiplicity $k$  gives rise to
logarithmic term up to $\log^{k-1}$. 

In particular this means that near any of the 25 roots of $Q_8(x)$ the local solutions
have no logarithmic terms and the solutions are thus {\em analytic} since all the
exponents are positive integers. The roots of $Q_8(x)$ are thus {\em apparent singularities}
\cite{ForsythV4,Ince} of the Fuchsian equation (\ref{eq:TCfde}). There are methods for 
distinguishing real and apparent singularities (see, e.g,  \cite{ForsythV4} \S 45)
and in principle one should check that the roots  of $Q_8(x)$ satisfy the conditions for 
being apparent singularities. However, this theoretical method is quite cumbersome. An easier 
numerical way to see that the roots of $Q_8(x)$ must be apparent singularities is as follows. 
We already found a 10th order Fuchsian equation for which the polynomial $P_{10}(x)$ was of a 
form similar to $P_8(x)$ as listed in equation (\ref{eq:TCpol}), but with the degree of $Q_{10}(x)$ 
being only 7. That is all the singularities as tabulated in Table~\ref{tab:TCexp} also appear
in this higher order equation with the exception of the 25 roots of $Q_8(x)$ (at most
7 of these could appear in the order 10 Fuchsian equation). In fact we can find a
solution of order 14 of the same form as above but with $Q_{14}(x)$ being just a constant.
So at this order none of the roots of $Q_8(x)$ appear. Clearly 
any real singularity
of the system cannot be made to vanish and we conclude that the  25 roots of $Q_8(x)$
must indeed be apparent singularities.

Assuming that only repeated roots give rise to logarithmic terms, and thus that a sequence
of positive integers give rise to {\em analytic} terms, then
near the physical critical point $x=x_c=1/4$ we expect the singular behaviour
\begin{equation}\label{eq:xc}
\TC(x) \sim A(x) (1-4x)^{-1/2} + B(x) (1-4x)^{-1/2}\log(1-4x),
\end{equation}
where $A(x)$ and $B(x)$ are analytic in the neighbourhood of $x_c$. Note that the terms
associated with the exponents $1/2$ and $3/2$ become part of the analytic correction to the
$(1-4x)^{-1/2}$ term.
Near the singularity on the negative $x$-axis, $x=x_-=-1/4$ we expect the singular behaviour
\begin{equation}\label{eq:xm}
\TC(x) \sim C(x) (1+4x)^{13/2},
\end{equation}
where again $C(x)$ is analytic near $x_-$. We expect similar behaviour near the pair of 
singularities $x=\pm \rmi/2$, and finally at the
roots of $1+x+7x^2$ we expect the behaviour $\TC(x) \sim D(x)(1+x+7x^2)^2\log (1+x+7x^2)$.

Next we turn our attention to the asymptotic behaviour of the coefficients of $\TC (x)$.
To standardise our analysis, we assume that the critical point is at 1. The
growth constant of staircase and imperfect staircase polygons is 4, so we normalise
the series by considering a new series with coefficients $r_n$, defined by 
$r_n = t_{n+2}/4^n.$
Thus the generating function we study is
$\RTC(y) = \sum_{n\geq 0} r_ny^n = 4 + 3y + 2.625y^2 + \cdots$.
From equations (\ref{eq:xc}) and (\ref{eq:xm}) 
it follows that the asymptotic form of the coefficients is
\begin{equation}\label{eq:asymp}
[y^n]\RTC(y) = r_n = \frac{1}{\sqrt{n}}\sum_{i \ge 0} \left(\frac{a_i\log{n} + b_i}{n^i}
+ (-1)^n\left( \frac{c_i}{n^{7+i}} \right) \right ) + {\rm O}(\lambda^{-n}).
\end{equation}
The last term includes the effect of other singularities, further from the origin than
the dominant singularities. These will decay exponentially since $\lambda > 1$ in the scaled variable $y=x/4$.

Using the recurrence relations for $t_n$ (derived from the ODE) it is easy and fast
to generate many more terms $r_n$. We generated the first 100000 terms and saved 
them as floats with 500 digit accuracy (this calculation took less than 15 minutes).
With such a long series it is possible to obtain accurate numerical estimates of the first 20 amplitudes 
$a_i$, $b_i$, $c_i$ for $i \le 19$  with precision of more than 100 digits for the dominant
amplitudes, shrinking to 10--20 digits for the the case when $i = 18,$ or $19.$ 
In making these estimates we have ignored the exponentially decaying term, which is the last term
in eq.(\ref{eq:xm}).
In this way we
confirmed an earlier conjecture \cite{CGD97} that $a_0 = \frac{3\sqrt{3}}{\pi^{3/2}},$
where we have taken into account the different normalisation, as discussed in the introduction.
We also find that $b_0 = 3.173275384589898481765\ldots$ and $c_0 = \frac{-24}{\pi^{3/2}}$,
though we have not been able to identify $b_0$. However, we have successfully identified further
sub-dominant amplitudes, and find $a_1=\frac{-89}{8\sqrt{3}\pi^{3/2}},$ 
$a_2=\frac{1019}{384\sqrt{3}\pi^{3/2}},$ and
$a_3=\frac{-10484935}{248832\sqrt{3}\pi^{3/2}},$ and
$c_1 = \frac{225}{\pi^{3/2}},$ 
$c_2 = \frac{-16575}{16\pi^{3/2}},$ and $c_3 = \frac{389295}{128\pi^{3/2}}.$ 
It seems possible that the amplitudes $\pi^{3/2}\sqrt{3}a_i$ and $\pi^{3/2}c_i$ are rational.

Estimates for the amplitudes were obtained by fitting $r_n$ to the form given above using
an increasing number of amplitudes. 
`Experimentally' we find we need about the same total number of terms at $x_c$ and $-x_c=x_-$.

So in the fits we used the terms with amplitudes $a_i$, and $b_i$, $i=0,\ldots,K$ and 
$c_i$, $i=0,\ldots,2K$. Going only to $i=K$ with the $c_i$ amplitudes results in much poorer
convergence and going beyond $2K$ leads to no improvement.
For a given $K$ we thus have to estimate $4K+3$ unknown amplitudes. So we use the
last $4K+3$ terms $r_n$ with $n$ ranging from 100000 to $100000-4K-2$ and solve the
resulting $4K+3$ system of linear equations.
We find that the amplitudes are fairly stable up to around $2K/3$.
We observed this by doing the calculation with $K=30$ and $K=40$ and then looking
at the difference in the amplitude estimates. For $a_0$ and $b_0$ the
difference is less than $10^{-131}$, while for $c_0$ the difference is less than $10^{-123}$.
Each time we increase the amplitude index by 1 we loose around $10^6$ in accuracy.
With $i=20$ the differences are respectively around $10^{-16}$ and $10^{-8}$.

The excellent convergence is solid evidence (though naturally not a proof) that
the assumptions leading to equation~(\ref{eq:asymp}) are correct. Further evidence
was obtained as follows. We can add extra terms to the asymptotic form and
check what happens to the amplitudes of the new terms. If the amplitudes are
very small it is highly likely that the terms are not truly present (if the
calculation could be done exactly these amplitudes would be zero). One possibility
is that our assumption about integer exponents leading only to analytic terms is incorrect.
To test this we fitted to the form
$$\frac{1}{\sqrt{n}}\sum_{i \ge 0} \left(\frac{\tilde{a}_i\log{n} + \tilde{b}_i}{n^{i/2}}
+ (-1)^n\left( \frac{\tilde{c}_i}{n^{7+i}} \right) \right ) + {\rm O}(\lambda^{-n}),
$$
(as above, in making these estimates we have ignored the exponentially decaying term, 
which is the last term in the above equation.)
With $K=30$ we found that the amplitudes $\tilde{a_1}$ and $\tilde{b_1}$ of the
terms $\log{n}/n$ and $1/n$, respectively, were less than $10^{-60}$, while the
amplitudes $\tilde{a_3}$ and $\tilde{b_3}$ were less than $10^{-50}$. We think
we can safely say that all the additional terms we just added are {\em not} present.
We found similar results if we added terms like $\log^2{n}$ or additional
$\log{n}$ terms at $y=-1$. That is, we found that those terms were not present.
 So this fitting procedure provides convincing
evidence that the asymptotic form (\ref{eq:asymp}), and thus the assumption
leading to this formula, is correct.

\section{Conclusion}

We have developed an improved algorithm for enumerating three choice polygons. The
extended series, coupled with a search program that assumes the solution is a {\em Fuchsian}
ODE, enabled us to discover the underlying ODE, which was of $10^{th}$ order. We
did this without using more than 50 of the coefficients that we had generated. That
is to say, 50 of the known coefficients were unused, and so their value provided
a check on the solution found. This
leads us to believe that we have found the correct ODE, as it reproduces the known, unused
coefficients, though we have not provided a proof.
Further refinement allowed us to find an $8^{th}$ order ODE. 

A numerical technique we have developed specifically for such problems then allowed us to 
find accurate numerical estimates for the amplitudes of the first several terms in the
asymptotic form for the coefficients.

We have also initiated an investigation of the {\em area} generating function. We expect this
to involve $q$-series, and thus far our investigations only lead us to believe that the area
generating function $A(q)$ is of the form
$$A(q) = (G(q) + H(q)/\sqrt{1 - q/\eta})/[J_0(1,1,q)^2],$$
where $J_0(x,y,q)$ is a $q$-generalisation of the Bessel function, and occurs, for
example, in the solution of the problem of staircase polygons enumerated by area \cite{BM96a}.
Here $q=\eta$ is the first zero of $J_0(1,1,q)$, and $G$ and $H$ are regular in the
neighbourhood of $q = \eta.$ The coefficients thus behave asymptotically as 
$$a_n = [q^n]A(q) \sim const. \eta^{-n}n^{3/2}.$$
In a subsequent publication \cite{GJ06} we propose to investigate the area generating function more
fully, and to say more about the properties of the ODE we have found for
the perimeter generating function. In particular, we discuss some simple solutions of the ODE, and 
ask what these can tell us about the full solution.

\section*{Acknowledgments}

We would like to thank N Zenine and J-M Maillard for access to their Mathematica
routines for identifying  differential equations and useful advice about their use.
We would also like to thank M. Bousquet-M\'elou for communicating her unpublished
results on this problem, which we mentioned in the introduction.
We gratefully acknowledge financial support from the Australian Research Council.

\section*{\label{app:TCpol} Appendix A: Polynomials $\bm Q_n(x)$}

\begin{eqnarray*}
 Q_8(x) &=& -180 + 4005x - 45340 x^2 + 352567 x^3 
   - 2100653 x^4 + 8247059 x^5 + 1869782 x^6  \\
  &&- 198745492 x^7 + 222232422 x^8 + 7981490552 x^9 - 58454247760 x^{10} \\ 
  &&+ 223070561538 x^{11} - 653903984242 x^{12} + 1691567153918 x^{13}  \\
  &&- 3628069390936 x^{14} + 9508812403200 x^{15} - 42130737708796 x^{16} \\
  &&+ 151950842991736 x^{17} - 347187650580720 x^{18} + 558723092175488 x^{19} \\
  &&- 722483977609792 x^{20} + 551434913787008 x^{21} + 82817126361856 x^{22} \\
  &&- 426478334005248 x^{23} + 279157576126464 x^{24} + 2780644737024 x^{25}
\end{eqnarray*}

\begin{eqnarray*}
 Q_7(x) &=& 
  -3420 + 82530 x - 926615 x^2 + 6866662 x^3 - 37878392 x^4 + 131975108 x^5 \\
&&  + 198512462 x^6 - 5322566116 x^7 + 16816064102 x^8 + 88956629348 x^9 \\
&&  - 872972184658 x^{10} + 3395585125316 x^{11} - 8662194926872 x^{12}  \\
&&  +  2179593948608 x^{13} + 130585482759744 x^{14} - 698610495175368x^{15} \\
&&  + 2229946022661696 x^{16} - 6216128747042864 x^{17} \\
&&  + 15724091332879132 x^{18} - 38607908490402392 x^{19} \\
&&  + 128963713249678592 x^{20} - 464640056155209952 x^{21} \\
&&  + 1296873363475699328 x^{22} -2966555758830491904 x^{23} \\
&&  + 5741739615453110784 x^{24}  -7824348079140616704 x^{25} \\
&&  + 8096625038421797888 x^{26} - 6327622359115208704 x^{27} \\
&&  - 663175049190105088 x^{28}  + 4390942020748738560 x^{29}\\
&&  - 3449431865352388608 x^{30} - 33011814317948928 x^{31}
\end{eqnarray*}

\begin{eqnarray*}
 Q_6(x) &=&
  -9180 + 310275 x - 4493475 x^2 + 40204094 x^3 - 262917778 x^4 \\
&&  + 1302960911 x^5 - 3743237840 x^6 - 8573351756 x^7 + 140454430666 x^8 \\
&&  - 409322626730 x^9 - 1570504457342 x^{10} + 18303308342032 x^{11} \\
&&  - 89658228463172 x^{12} + 259420736216632 x^{13} + 26862202296376 x^{14}\\ 
&&  - 4190021721023184 x^{15} + 21897720821926584 x^{16} \\
&&  - 75837533674259508 x^{17} + 212508813586272428 x^{18} \\
&&  - 476010656497826944 x^{19} + 1034090705056496672 x^{20} \\
&&  - 3196181326637410304 x^{21} + 10833216991064882848 x^{22} \\
&&  - 30172750280212408832 x^{23} + 70340668591569812736 x^{24} \\
&&  - 132300506186507025408 x^{25} + 177280513560453634560 x^{26} \\
&&  - 184990945124657242112 x^{27} + 135828858351882342400 x^{28} \\
&&  + 12754320650381836288 x^{29} - 85576383794502107136 x^{30} \\
&&  + 61165017902554546176 x^{31} + 565560894352785408 x^{32}
\end{eqnarray*}

\begin{eqnarray*}
 Q_5(x) &=&-4500 + 244800 x - 4876845 x^2 + 55164150 x^3 - 438701640 x^4 \\
&&  + 2758453094 x^5 - 13804842198 x^6 + 45370091528 x^7 - 3608230380 x^8 \\
&&  - 892524490064 x^9 + 4421327158154 x^{10} - 2297315126532 x^{11} \\
&&  - 103201897035096 x^{12} + 748998082407080 x^{13} - 2329708885595260 x^{14} \\
&&  - 457382726191024 x^{15} + 35817660448173240 x^{16} \\
&&  - 188156345496838984 x^{17} + 677783573996257364 x^{18} \\
&&  - 1904649390940935752 x^{19} + 4199594693024922016 x^{20} \\
&&  - 8814226144821806432 x^{21} + 23568486792872894272 x^{22} \\
&&  - 70089404940793421632 x^{23} + 188311273940137111552 x^{24}  \\
&&  - 435002993494719438848 x^{25} + 791152555777632593920 x^{26} \\
&&  - 1045593345640931730432 x^{27} + 1096015208846337957888 x^{28} \\
&&  - 774016903940080771072 x^{29} - 37178029375778357248 x^{30} \\
&&  + 412071049964952354816 x^{31} - 275345921075326746624 x^{32} \\
&&  - 2464051649845395456 x^{33}
\end{eqnarray*}

\begin{eqnarray*}
 Q_4(x) &=&31500 - 1114080 x + 17560755 x^2 - 178469565 x^3 + 1412918104 x^4 \\
&&  - 9431590849 x^5 + 52336335969 x^6 - 220707961458 x^7 + 525965711332 x^8 \\
&&  + 915935968370 x^9 - 13996439933349 x^{10} + 35303141246088 x^{11}  \\
&&  + 231992664240696 x^{12} - 2180352456480752 x^{13} + 6859298731027888 x^{14} \\
&&  - 1272161420555012 x^{15} - 75338205421491734 x^{16}  \\
&&  + 406836568590013948 x^{17} - 1513874477368697252 x^{18} \\
&&  + 4439738234446975124 x^{19} - 10514406278248398472 x^{20} \\
&&  + 22797584086521040520 x^{21} - 52624647215757093584 x^{22} \\
&&  + 130673617185226821792 x^{23} - 324680301683722155712 x^{24} \\
&&  + 724969297042825531136 x^{25} - 1271869215082051692800 x^{26} \\
&&  + 1661614177465373698560 x^{27} - 1744486537247742479360 x^{28} \\
&&  + 1209724637255295010816 x^{29} - 32142663152460406784 x^{30} \\
&&  - 498040622799430975488 x^{31} + 321317702703841148928 x^{32} \\
&&  + 2787318284392857600 x^{33}
\end{eqnarray*}

\begin{eqnarray*}
 Q_3(x) &= &-156000 + 3778920 x - 46727325 x^2 + 457371630 x^3 - 3919246431 x^4 \\
&& + 27446185200 x^5 - 152613919692 x^6 + 659637747242 x^7 - 1723470963068 x^8 \\
&& - 1667066145852 x^9 + 27889854017778 x^{10} + 15933308039400 x^{11} \\
&& - 972460279627326 x^{12} + 4552136023731292 x^{13} - 7976188460233924 x^{14} \\
&& - 4422880527966948 x^{15} + 63325989574562728 x^{16} \\
&& - 287206984863975352 x^{17} + 1115308575007981980 x^{18} \\
&& - 3508943115779966584 x^{19} + 8987842561562515768 x^{20} \\
&& - 19184807012355087408 x^{21} + 37821550927408731776 x^{22} \\
&& - 83609060792238083072 x^{23} + 194683017390969665280 x^{24} \\
&& - 366258975512082319872 x^{25} + 498254429378056694784 x^{26}  \\
&& - 558421919820222289920 x^{27} + 441211762632912959488 x^{28} \\
&& - 80404063142199537664 x^{29} - 110342796490113417216 x^{30} \\
&& + 85904442856027127808 x^{31} + 720965567415582720 x^{32}
\end{eqnarray*}

\begin{eqnarray*}
 Q_2(x) &= & 102000 - 1245240 x + 445275 x^2 + 77507430 x^3 - 505005638 x^4 + 674357270 x^5 \\
&& + 7410398802 x^6 - 50751541108 x^7 + 109730141494 x^8 + 263567061768 x^9 \\
&& - 2398666258514 x^{10} + 4447124418524 x^{11} + 33544348232760 x^{12} \\
&& - 341405641395740 x^{13} + 1843130781900080 x^{14} - 7441271357292384 x^{15} \\
&& + 23827305830694324 x^{16} - 59142500096057112  x^{17} \\
&& + 113845825936073424 x^{18} - 169659492965796928 x^{19} \\
&& + 190085091157739584 x^{20} - 160391840217609984 x^{21} \\
&& + 95477320250924800 x^{22} - 21461546279272960 x^{23} \\
&& - 73590898428536832 x^{24} + 43442402559821824 x^{25} \\
&& + 129164030193680384 x^{26} - 136460131311329280 x^{27} \\
&& + 54532752690511872 x^{28} + 389290263183360 x^{29}
\end{eqnarray*}

\begin{eqnarray*}
 Q_0(x) =  Q_1(x)=-Q_2(x)
\end{eqnarray*}

\end{document}